\documentclass{elsart}

\usepackage{amsfonts}
\usepackage{amstext}
\usepackage{amssymb,amsmath}
\usepackage{graphicx}

\def\R{\mathbb{R}}

\def\0{\mathbf{0}}
\def\1{\mathbf{1}}
\def\cc{\mathbf{c}}

\def\r{\mathbf{r}}
\def\x{\mathbf{x}}
\def\X{\mathbf{X}}
\def\uu{\mathbf{u}}
\def\v{\mathbf{v}}
\def\w{\mathbf{w}}
\def\y{\mathbf{y}}
\def\z{\mathbf{z}}
\def\Y{\mathbf{Y}}
\def\eps{\varepsilon}

\def\I{\mathbf{I}}

\def\W{\mathbf{W}}
\def\A{\mathbf{A}}

\def\CC{\mathbf{C}}
\def\DD{\mathbf{D}}
\def\EE{\mathbf{E}}

\def\tr{\mathtt{tr}\,}
\def\diag{\mathtt{diag}\,}
\def\Span{\mathtt{Span}\,}

\def\Vol{{\mathtt{Vol}}}
\def\dist{{\mathtt{dist}}}

\def\part{\cal P} 
 
\def\tu{\tilde \uu}
\def\tv{\tilde \v}
\def\hu{\hat \uu}
\def\hv{\hat \v}

\begin{document}

\begin{frontmatter}
\title {SVD, discrepancy, and regular structure of contingency tables}
\author{Marianna Bolla\corauthref{ca}}
\ead{marib@math.bme.hu}
\corauth[ca]{Research supported in part by the Hungarian 
National Research Grants OTKA 76481 and OTKA-KTIA 77778; further, by the 
T\'AMOP-4.2.2.B-10/1-2010-0009 project.}

\address{Institute of Mathematics,
Budapest University of Technology and Economics}

\begin{abstract} 
We will use the factors obtained by correspondence analysis to find biclustering
of a contingency table such that the row--column cluster pairs are regular,
i.e., they have small discrepancy. In our main theorem,
the constant of the so-called volume-regularity is related to the SVD
of the normalized contingency table. Our result is applicable to two-way
cuts when both the rows and columns are divided into the same number of
clusters, thus extending partly the result of~\cite{Butler} estimating the 
discrepancy of a contingency table
by the second largest singular value of the normalized table (one-cluster,
rectangular case),
and partly the result of~\cite{Bolla11} for estimating the constant of 
volume-regularity by the structural eigenvalues and the distances of
the corresponding eigen-subspaces of the
normalized modularity matrix of an edge-weighted graph 
(several clusters, symmetric case). 
\end{abstract}

\begin{keyword}
Normalized contingency table
\sep Regular row-column pairs
\sep Biclustering
\sep Discrepancy
\sep Cluster variances
\sep Directed graphs

\end{keyword}

\end{frontmatter}

\section{Introduction}

A typical problem of contemporary cluster analysis  is to find 
relatively small number of groups of objects, belonging to
rows and columns of a contingency table which exhibit homogeneous
behavior with respect to each other and
do not differ significantly in size.
To make inferences on the separation that can be achieved for a given number
of clusters, minimum normalized two-way cuts are investigated and
related to the SVD of the correspondence matrix.

Contingency tables are
rectangular arrays with nonnegative, real entries. One example is 
the keyword--document matrix. Here the entries are associations between
documents and words. Based on network data, the entry in the $i$th row
and $j$th column is the relative frequency of word $j$
in document $i$. Latent semantic indexing looks
for real scores of the documents and keywords such that the score of a
any document be proportional to the total scores of the keywords occurring
in it, and vice versa, the score of any keyword be proportional to the total
scores of the documents containing it. Not surprisingly, the solution is given 
by the SVD of the binary table, where the document- and keyword-scores are the
coordinates of the left and right singular vectors corresponding to its largest
non-trivial singular value which gives the constant of proportionality. 

This idea is generalized in~\cite{FKV} in the following way.
We can think of the above relation between keywords and documents
as the relation with respect to the most important topic (or context, or
factor). After this, we are looking for another scoring with respect to
the second topic, up to $k$ (where $k$ is a positive integer not exceeding the 
rank of the table). The solution is given by the singular vector pairs
corresponding to the $k$ largest singular values of the table.

If a scoring system
is endowed with the marginal measures, the problem can be formulated in
terms of correspondence analysis and correlation maximization. The
problem is solved by the SVD of the correspondence matrix (normalized
contingency table), where the singular vector pairs are also transformed,
see~\cite{Bol6}.
In this way, instead of scores,
the documents and keywords have $k$-dimensional representatives, based of which
further investigations, spacial representation, or biclustering can be 
performed that finds simultaneous clustering of the rows and columns 
of the table with densities as 
homogeneous as possible between the keyword--document cluster pairs.

The problem is also related to the Pagerank (see~\cite{Kleinberg}) and
to microarray analysis (see~\cite{Kluger}) 
when we want to find clusters of the rows and columns of
a microarray, simultaneously. Here rows
correspond to genes and columns to different conditions, whereas the
entries are expression levels of genes under specific
conditions.
We also look for a bipartition of the genes and conditions
such that genes in the same cluster equally (not necessarily weakly
or strongly) influence conditions of the same cluster.

In Section~\ref{svdcont} we deal with the singular value decomposition (SVD) 
of a correspondence matrix.
In Section~\ref{norm} we relate it to normalized two-way cuts of the 
contingency
table, while in Section~\ref{reg} the constant of volume-regularity of
row--column clusters pairs is estimated by means of the SVD.
Section~\ref{discussion} is devoted to discussion, application and possible
extension to directed graphs.

\section{ SVD of contingency tables and correspondence matrices}\label{svdcont}

 Let $\CC$
be a contingency table on row set $Row =\{ 1,\dots ,n \}$ and column set 
$Col =\{ 1,\dots ,m \}$, where $\CC$ is $n\times m$
matrix of entries $c_{ij} \ge 0$.
Without loss of generality, we suppose that
there are not identically zero rows or columns.
 Here $c_{ij}$ is some kind of association
between the objects behind row $i$ and column $j$, where 0 means no 
interaction at all. 

Let the row- and column-sums of $\CC$ be 
$$
 d_{row,i} =\sum_{j=1}^m c_{ij}  \quad (i=1,\dots ,n ) \quad \mathrm{and}\quad
 d_{col,j} =\sum_{i=1}^n c_{ij}  \quad (j=1,\dots ,m )
$$
which are
collected  in the main diagonals of the $n\times n$ and $m\times m$
diagonal matrices $\DD_{row}$ and $\DD_{col}$, respectively.

For a given integer $1\le k\le \min \{ n,m \}$,
 we are looking for $k$-dimensional
representatives $\r_1 ,\dots ,\r_n$ of the rows and 
 $\cc_1 ,\dots ,\cc_m$ of the columns
such that they minimize the objective function
\begin{equation}\label{objrec}
 Q_k = \sum_{i=1}^n \sum_{j=1}^m  c_{ij} \| \r_i -\cc_j \|^2 
\end{equation}
subject to
\begin{equation}\label{consts}
 \sum_{i=1}^n d_{row,i} \r_i \r_i^T =\I_k , \quad
 \sum_{j=1}^m d_{col,j} \cc_j \cc_j^T =\I_k .
\end{equation}

When minimized, the objective function $Q_k$ favors $k$-dimensional placement
of the rows and columns such that representatives of highly associated 
rows and columns are
forced to be close to each other. As we will see, this is equivalent to the
problem of correspondence analysis. 

Indeed, let us put both the objective function and the constraints in a more
favorable form. Let $\X$ be the $n\times k$ matrix of rows $\r_1^T ,\dots,
\r_n^T$; let $\x_1 ,\dots ,\x_k \in \R^n$ denote the columns of $\X$, 
for which fact we use the notation $\X =(\x_1 ,\dots ,\x_k )$.
Similarly, let $\Y$ be the $m\times k$ matrix of rows $\cc_1^T ,\dots,
\cc_m^T$; let $\y_1 ,\dots ,\y_k \in \R^m$ denote the columns of $\Y$, i.e.,
$\Y =(\y_1 ,\dots ,\y_k )$. 
Hence, the
constraints~(\ref{consts}) can be formulated  like 
$$
 \X^T \DD_{row} \X =\I_k, \quad \Y^T \DD_{col} \Y =\I_k .
$$

With this notation, the objective function~(\ref{objrec}) is
\begin{equation}\label{Qrec}
\begin{aligned}
 Q_k &=  \sum_{i=1}^n \sum_{j=1}^m  c_{ij} \| \r_i -\cc_j \|^2  =
 \sum_{i=1}^n d_{row,i} \| \r_i \|^2 +\sum_{j=1}^m d_{col,j} \| \cc_j \|^2 -
\sum_{i=1}^n \sum_{j=1}^m  c_{ij} \r_i^T \cc_j    \\
& = 2k- \tr \X^T \CC \Y = 2k- \tr (\DD_{row}^{1/2} \X )^T 
 (\DD_{row}^{-1/2} \CC  \DD_{col}^{-1/2} ) 
  (\DD_{col}^{1/2} \Y ) ,
\end{aligned}
\end{equation}
where the matrix $\CC_{corr} =\DD_{row}^{-1/2} \CC  \DD_{col}^{-1/2} $ is the
\textit{correspondence matrix (normalized contingency table)} 
belonging to the table $\CC$,
see~\cite{Bol6}. 
If we multiply all the entries of $\CC$ with the same positive constant,
the correspondence matrix $\CC_{corr}$ will not change. Therefore, without
the loss of generality, $\sum_{i=1}^n \sum_{j=1}^m c_{ij}=1$ will be supposed
in the sequel.
The correspondence matrix has SVD
\begin{equation}\label{svd}
 \CC_{corr} =\sum_{i=1}^r s_i \v_i \uu_i^T ,
\end{equation}
where $r\le \min \{ n,m \}$ is the rank of $\CC_{corr}$, or equivalently
(since there are not identically zero rows or columns), the rank of $\CC$.
Here $1=s_1 \ge s_2 \ge \dots \ge s_r >0$
are the non-zero singular values of $\CC_{corr}$, and 1 is a single singular 
value if $\CC_{corr}$, or equivalently,  $\CC$ is non-decomposable
($\CC \CC^T$ is irreducible).
In this case $\v_1 =(\sqrt{d_{row,1} }, \dots ,\sqrt{d_{row,n} } )^T$ and
$\uu_1 =(\sqrt{d_{col,1} }, \dots ,\sqrt{d_{col,m} } )^T$.

Note that the singular spectrum of a decomposable contingency table
can be composed from the singular spectra of its non-decomposable parts, as
well as their singular vector pairs. Therefore, in the future, the 
non-decomposability of the underlying contingency table will be supposed.
In this way, the following representation theorem for contingency tables
can be formulated.

\begin{thm}\label{repr}
Let  $\CC$ be a non-decomposable contingency table 
with SVD~(\ref{svd}) of its  correspondence matrix $\CC_{corr}$. 
Let $k \le r$ be a positive integer such that 
$s_k > s_{k+1}$. Then the minimum of~(\ref{objrec}) subject to~(\ref{consts}) 
is $2k -\sum_{i=1}^k s_i$ and it is attained with
the optimum row representatives $\r_1^* ,\dots ,\r_n^*$ 
and column representatives $\cc_1^* ,\dots ,\cc_m^*$, 
the transposes of which are
row vectors of $ \X^* = \DD_{row}^{-1/2} (\v_1 ,\dots ,\v_k )$
and $\Y^* =\DD_{col}^{-1/2} (\uu_1 ,\dots ,\uu_k )$, respectively.
\end{thm}

\begin{pf}
In view of~(\ref{Qrec}), we have to maximize
$$
 \tr (\DD_{row}^{1/2} \X )^T \CC_{corr}  (\DD_{col}^{1/2} \Y )
$$
under the given constraints. Separation theorems for the singular value
decomposition (see e.g.,~\cite{Bhatia} and~\cite{Rao}) 
are applicable, yielding the required statement.
\end{pf}

The vectors $\r_1^* ,\dots ,\r_n^*$ and $\cc_1^* ,\dots ,\cc_m^*$ giving the 
optimum in the above theorem are
called \textit{optimum k-dimensional representatives} of the rows and columns, 
while the
transformed singular vectors 
$\DD_{row}^{-1/2} \v_1 ,\dots ,\DD_{row}^{-1/2} \v_k$ and
$\DD_{col}^{-1/2} \uu_1 ,\dots ,$ $\DD_{col}^{-1/2} \uu_k $ are called 
\textit{vector components} of the rows and columns
taking part in the $k$-dimensional representation.

Observe that
the dimension $k$ does not play an important role here: the vector components
can be included successively up to a $k$ such that $s_k >s_{k+1}$.
We remark that the singular vectors can arbitrarily be chosen in the 
isotropic subspaces
corresponding to possible multiple singular values, under the orthogonality 
conditions.
Further, provided that 1 is  a single singular value, the first 
vector components are the constantly $\1$
vectors in $\R^n$ and $\R^m$, respectively, and hence, the $k$-dimensional
representation is realized in a $(k-1)$-dimensional hyperplane of $\R^k$.

A symmetric contingency table corresponds to a weighted graph, and our
correspondence matrix is the identity minus the normalized Laplacian,
called normalized modularity matrix in~\cite{Bolla11}.
In another view, a contingency table can be considered as part of the
weight matrix of a bipartite graph on vertex set $Row \cup Col$.
However, it would be tedious to always distinguish between these two types
of vertices, we rather use the framework of correspondence analysis,
and formulate our statements in terms of rows and columns.

\section{Normalized two-way cuts of contingency tables}\label{norm}

Given the $n\times m$ contingency table $\CC $ on row set $Row$ and column set
$Col$, further, an integer $k$ $(0<k\le r)$, 
we want to simultaneously partition its rows and columns into  
disjoint, nonempty subsets 
$$
 Row =R_1 \cup \dots \cup R_k , \quad Col =C_1 \cup \dots \cup C_k
$$
such that the \textit{cut}s $c(R_a ,C_b ) =\sum_{i\in R_a} \sum_{j\in C_b} c_{ij}$
$(a,b=1,\dots ,k)$ between the row-column cluster pairs
be as homogeneous as possible. For this requirement, the following so-called 
\textit{normalized two-way cut} of the contingency table with respect to 
the above $k$-partitions $P_{row} =(R_1 ,\dots ,R_k )$ and $P_{col} =(C_1 ,\dots
,C_k )$ of its rows and columns and the collection of signs $\sigma$ is defined
as follows:
$$
 \nu_{k} (P_{row} ,P_{col} ,\sigma )=  
\sum_{a=1}^k \sum_{b=1}^k \left( \frac1{\Vol (R_a )} +
    \frac1{\Vol (C_b )} +
 \frac{2\sigma_{ab} \delta_{ab} }
 {\sqrt{\Vol (R_a )\Vol (C_b )}}  \right)  c(R_a ,C_b ) ,
$$
where 
$$
 \Vol (R_a ) =\sum_{i\in R_a} d_{row,i}=\sum_{i\in R_a} \sum_{j=1}^m c_{ij} ,\quad
 \Vol(C_b ) = \sum_{j\in C_b} d_{col,j}=\sum_{j\in C_b} \sum_{i=1}^n c_{ij}
$$
are volumes of the clusters, $\delta_{ab}$ is the Kronecker delta, and 
the sign $\sigma_{ab}$ is equal to 1 or -1 (it only has relevance in the
$a=b$ case, when it helps balancing between the volumes of the same index
row and column clusters), $\sigma:= (\sigma_{11} ,\dots ,\sigma_{kk})$.
We want to minimize the above normalized two-way cut  with respect to
all possible $k$-partitions ${\part}_{row,k}$ and
${\part}_{col,k}$ of the rows and columns, further, to $\sigma$, simultaneously.
The objective function penalizes row- and column clusters of extremely 
different volumes in the $a\ne b$ case, whereas in the $a=b$ case $\sigma_{aa}$
moderates the balance between $\Vol (R_a )$ and $\Vol (C_a )$.

\begin{defn}
The normalized two-way cut of the contingency table $\CC $ is
$$
  \nu_k (\CC ) =\min_{P_{row} ,P_{col}, \sigma}\nu_{k}(P_{row},P_{col},\sigma ).
$$
\end{defn}

\begin{thm} Let $1=s_1 > s_2 \dots \ge s_r$ be the positive singular values of
the correspondence matrix belonging to the non-decomposable contingency table
$\CC$ of rank $r$, and $k\le r$ be a positive integer. Then  
$$
 \nu_k (\CC ) \ge 2k -\sum_{i=1}^k s_i .
$$
\end{thm}

\begin{pf}
We will show that $\nu_{k} (P_{row} ,P_{col} ,\sigma )$ is $Q_k$ in the special 
representation, where the column vectors of $\X$
and $\Y$ are partition vectors belonging to $P_{row}$ and $P_{col}$,
respectively. Therefore, the statement follows, as
the overall minimum is $2k -\sum_{i=1}^k s_i$. 
Indeed, let the $i$th coordinate of the left vector component
$\x_a$ be
$$
 x_{ia} := \frac1{\sqrt{\Vol (R_a )}}\quad \textrm{if} \quad  
       i\in R_a ,\ a=1,\dots k;
$$
similarly, let the $j$th coordinate of the right vector 
component $\y_b$ be
$$
 y_{jb} =\sigma_{bb} \frac1{\sqrt{\Vol (C_b )}} \quad \textrm{if} \quad  
        j\in C_b, \ b=1,\dots ,k,
$$
otherwise the coordinates are zeros.
With this, the matrices $\X$ and $\Y$ satisfy the conditions
imposed on the representatives, further
$$
 \| \r_i -\cc_j \|^2 = \frac1{\Vol (R_a )} +\frac1{\Vol (C_b )}
 +\frac{2\sigma_{bb} \delta_{ab}}  
 {\sqrt{\Vol (R_a )\Vol (C_b )}}  , \quad
 \text{if} \quad i\in R_a , \, j\in C_b. 
$$
\end{pf}

In case of a symmetric contingency table (weight matrix $\W$ of an 
edge-weighted graph), we get the same
result with the representation based on the eigenvectors 
belonging to the largest absolute value eigenvalues of the 
matrix $\DD^{-1/2} \W \DD^{-1/2}$, where $\DD =\DD_{row}=\DD_{col}$, 
see~\cite{Bolla11}. 
However, $\nu_{k} (P_{row} ,P_{col} ,\sigma )$ cannot always
be directly related to the normalized cut, except the following
two special cases.
\begin{itemize}
\item
When the $k-1$ largest absolute value
eigenvalues of the normalized modularity matrix are all positive, or
equivalently, if the $k$ smallest eigenvalues (including the zero) of the
normalized Laplacian matrix are farther from 1 than any other  eigenvalue
which is greater than 1. In this case the $k-1$ largest singular values 
(apart from the 1) of the correspondence matrix are identical to the $k-1$
largest eigenvalues of the normalized modularity matrix, and the left and right
singular vectors are identical to the corresponding eigenvector with the
same orientation. Consequently, for the $k$-dimensional
row- and column-representatives $\r_i =\cc_i$
($i=1,\dots ,n=m$) holds.
With the choice $\sigma_{bb} =1$ $(b=1,\dots ,k )$, 
the corresponding $\nu_k (\CC )$ is twice the 
normalized cut of
our weighted graph in which weights of edges within the clusters do not count.
In this special situation, the normalized two-way cut
also favors $k$-partitions with low inter-cluster edge-densities
(consequently, intra-cluster densities tend to be large, as they do not
count in the objective function).
\item
When the $k-1$ largest absolute value eigenvalues of the normalized
modularity matrix are all negative, 
then $\r_i =-\cc_i$ for all $(k-1)$-dimensional
row and column representatives, and any (but only one)
of them can be the corresponding vertex representative. Now $\nu_k (\CC )$,
which is attained with the choice  $\sigma_{bb} =-1$ $(b=1,\dots ,k )$,
differs from the normalized cut in that it also counts the edge-weights
within the clusters. 
Indeed, in the $a=b$, $R_a =C_a =V_a$ case
$$
 \| \r_i -\cc_j \|^2 = \frac1{\Vol (V_a )} +\frac1{\Vol (V_b )}   
 +\frac2{\sqrt{\Vol (V_a )\Vol (V_b )}}  = 
 \frac4{\Vol (V_a )} 
$$
if $i,j\in V_a$. 
Here, by minimizing the normalized $k$-way cut,  rather a so-called 
anti-community structure is detected in that $c(R_a ,C_a )=c(V_a ,V_a )$
is suppressed to compensate for the term $\frac4{\Vol (V_a )}$.
\end{itemize}

We remark that Ding et al.~\cite{Ding} treat this problem for two row- and
column-clusters and minimize another objective function such that it
favors 2-partitions where $c (R_1 ,C_2 )$ and $c (R_2 ,C_1 )$ are small
compared to $c (R_1 ,C_1 )$ and $c (R_2 ,C_2 )$. The solution is also given
by the transformed $\v_2 , \uu_2$ pair. However, it is the objective
function $Q_k$ which best complies with the SVD of the correspondence matrix, 
and hence, gives the continuous relaxation of the normalized cut minimization
problem. The idea of Ding et al.
could be naturally extended to the case of several, but the same number of
 row and column
clusters, and it may work well in the keyword-document classification problem.
Though, in some real-life problems, e.g., clustering genes and 
conditions
of microarrays, we rather want to find clusters of similarly functioning genes
that equally (not especially weakly or strongly) influence conditions of the 
same cluster.
Dhillon~\cite{Dhillon} also suggests a multipartition algorithm that runs
the k-means algorithm simultaneously for the row and column representatives.

\section{Regular row-column cluster pairs}\label{reg}

Let us start with the one-cluster case. 
Let $\CC $ be an $n\times m$ contingency table and let $\CC_{corr}$ be
the correspondence matrix 
belonging to it. The Expander Mixing Lemma for edge-weighted graphs 
naturally extends to this situation, see the following result of~\cite{Butler}.

\begin{prop}\label{expcont}
Let $\CC $ be a non-decomposable contingency table (i.e., $\CC \CC^T$ is  
irreducible) on row set $Row$ and column 
set $Col$, and of  total volume 1. 
Then for all $R\subset Row$ and $C\subset Col$
$$
 | c (R, C) - \Vol (R) \Vol (C)| 
 \le s_2  \sqrt{\Vol (R) \Vol (C)} ,
$$
where $s_2$ is the largest but 1 singular value of the normalized contingency 
table $\CC_{corr}$. 
\end{prop}

Since the spectral gap of $\CC_{corr}$ is  $1-s_2$, in view of the above 
Expander Mixing Lemma,  'large' spectral gap
is an indication that the weighted cut between any 
row and column subset of the contingency table is near to 
what is expected in a random table.
The following notion of discrepancy is just measures the deviation from
this random situation.  
The discrepancy (see~\cite{Butler}) of the contingency table  $\CC$ of
total volume 1 is the smallest $\alpha >0$ such that for all $R \subset Row$
and $C\subset Col$
$$
 | c (R, C) - \Vol (R) \Vol (C)| \le \alpha  \sqrt{\Vol (R) \Vol (C)} .
$$
In view of this, the result of Theorem~\ref{expcont} can be interpreted as
follows:
$\alpha$ singular value separation causes $\alpha$ discrepancy, where the 
singular value
separation is the second largest singular value of the normalized 
contingency table, 
which is the smaller the bigger the separation between the 
largest singular value (the 1) of the normalized contingency table 
and the other singular values is. Based on the ideas  of 
\cite{Bilu} and \cite{Bollobas}, Butler~\cite{Butler} proves the converse of 
the Expander Mixing Lemma for contingency tables, namely that
$$
 s_2 \le 150\alpha (1-8\log\alpha ) .
$$
Now we extend the notion of discrepancy to volume-regular pairs. 
\begin{defn}\label{volregcont}
The row--column cluster pair $R\subset Row$, $C\subset Col$ of the 
contingency table $\CC$ of total volume 1 
is $\gamma$-volume regular
if for all $X\subset R$ and $Y\subset C$ the relation
\begin{equation}\label{jeles}
| c (X, Y) -\rho (R,C) \Vol (X) \Vol (Y)| \le \gamma \sqrt{\Vol (R) \Vol (C)} 
\end{equation}
holds, where $\rho (R,C) =\frac{c(R,C)}{ \Vol (R) \Vol (C)}$ is the relative
inter-cluster density of the row--column pair $R,C$.
\end{defn}

Now we will show that for given $k$, if the clusters are formed via
applying the weighted $k$-means algorithm for the optimal row- and column 
representatives,
respectively, then the so obtained row--column cluster pairs are homogeneous
in the sense that they form equally dense parts of the contingency table.
More precisely, the constant $\gamma$ of the volume regularity of the
pairs will be related to the SVD of $\CC_{corr}$.
To this end, we introduce the following notion.

The weighted $k$-variance of the $k$-dimensional row representatives 
is defined by
\begin{equation}\label{kszoras}
 {S}_k^2 (\X ) =
\min_{(R_1 ,\dots ,R_k )}
\sum_{a=1}^k \sum_{j\in R_a } d_{row,j} \| \r_j -{\bar \r}_a \|^2 ,
\end{equation}
where ${\bar \r}_a =\frac1{\Vol (R_a ) } \sum_{j\in R_a } d_{row,j} \r_j $ is the
weighted center of cluster $R_a$  $(a=1,\dots ,k )$. Similarly,
the weighted $k$-variance of the $k$-dimensional column representatives is 
\begin{equation}\label{kszoras1}
 {S}_k^2 (\Y ) =
\min_{(C_1 ,\dots ,C_k )}
\sum_{a=1}^k \sum_{j\in C_a } d_{col,j} \| \cc_j -{\bar \cc }_a \|^2 ,
\end{equation}
where ${\bar \cc}_a =\frac1{\Vol (C_a ) } \sum_{j\in C_a } d_{col,j}\cc_j $ is the
weighted center of cluster $C_a$  $(a=1,\dots ,k )$.
Observe, that the trivial vector components can be omitted, and the 
$k$-variance of the so obtained $(k-1)$-dimensional representatives will be
the same.

\begin{defn}
The cut-norm of the rectangular real matrix $\A$ with row-set $Row$ and 
column-set $Col$ is
$$
 \| \A \|_{\square } =\max_{R\subset Row, \, C\subset Col}
  \left| \sum_{i\in R} \sum_{j\in C} a_{ij} \right| .
$$
\end{defn}

\begin{lem}\label{lemmam}
For the cut-norm of the $n\times m$ real matrix $\A$ 
$$
\| \A \|_{\square } \le \sqrt{nm} \| \A \|  
$$
holds, where the right hand side contains its spectral norm, i.e., the largest 
singular value of $\A$.
\end{lem}

\begin{pf}
$$
\aligned
 \| \A \|_{\square } &=\max_{\x \in \{ 0,1 \}^n,\, \y \in \{ 0,1 \}^m} 
 | \x^T \A \x | = 
\max_{\x \in \{ 0,1 \}^n ,\, \y \in \{ 0,1 \}^m} |
 (\frac{\x}{\| \x \|})^T \A 
(\frac{\x}{\| \x \|}) | \cdot \| \x \| \cdot \|\y\| | \\
&\le \sqrt{nm}
\max_{\| \x \| =1, \, \|\y\|=1} | \x^T \A \x | =\sqrt{nm} \| \A \| ,
\endaligned
$$
since for $\x \in \{ 0,1 \}^n$, $\| \x \| \le\sqrt{n}$, and
for $\y \in \{ 0,1 \}^m$, $\| \y \| \le\sqrt{m}$.
\end{pf}

\begin{thm}\label{tetelem}
Let $\CC$ be a non-decomposable contingency table of $n$-element row set $Row$ 
and $m$-element column set $Col$,  with row- and column sums
$d_{row,1} ,\dots ,d_{row,n}$ and $d_{col,1} ,\dots ,d_{col,m}$, respectively. 
Suppose that
$\sum_{i=1}^n \sum_{j=1}^m c_{ij} =1$ and there are no dominant rows and 
columns: $d_{row,i} =\Theta (1/n )$, $(i=1,\dots ,n)$ and
$d_{col,j} =\Theta (1/m )$, $(j=1,\dots ,m)$ as $n,m\to\infty$.   
Let the singular values  of $\CC_{corr}$ be 
$$
 1=s_1 >s_2 \ge \dots \ge s_k  >\eps \ge s _i  , 
 \quad i\ge k+1 .
$$
The partition $(R_1,\dots ,R_k )$ of $Row$ and 
$(C_1,\dots ,C_k )$ of $Col$ are defined so that they minimize
the weighted k-variances $S_k^2 (\X )$ and $S_k^2 (\Y )$ of the 
row and column representatives
defined in~(\ref{kszoras}) and (\ref{kszoras1}), respectively.
Suppose that there are  constants 
$0<K_1 ,K_2\le \frac1{k}$ such that $|R_i |\ge K_1 n$ and 
$|C_i |\ge K_2 m$ $(i=1,\dots ,k)$, respectively.
Then the $R_i, C_j$ pairs are
${\cal O} (\sqrt{2k} (S_k (\X ) S_k (\Y )) +\eps)$-volume regular 
$(i,j=1,\dots ,k)$.
\end{thm}

\begin{pf}
Recall that provided $\CC$ is non-decomposable, 
the largest singular value $s_1=1$ 
of  $\CC_{corr}$ is single with corresponding singular vector pair
$\v_1 =\DD_{row}^{1/2} \1$ and $\uu_1 =\DD_{col}^{1/2} \1$ with the constantly
$\1$ vectors of appropriate size. 
The optimal $k$-dimensional representatives of the rows and columns are row 
vectors  
of the matrices $\X =(\x_1 ,\dots, \x_k )$ and 
$\Y =(\y_1 ,\dots, \y_k )$, where $\x_i =\DD_{row}^{-1/2} \v_i$ and
$\y_i =\DD_{col}^{-1/2} \uu_i$, respectively 
$(i=1,\dots ,k)$. 
Suppose that the minimum   $k$-variance is attained on the
$k$-partition $(R_1 ,\dots ,R_k )$ of the rows and
$(C_1 ,\dots ,C_k )$ of the columns.
By an easy analysis of variance argument of~\cite{Bolla} it follows that
$$
 S_k^2 (\X ) =\sum_{i=1}^k \dist^2 (\v_i , F ) , \quad
 S_k^2 (\Y ) =\sum_{i=1}^k \dist^2 (\uu_i , G ) ,
$$
where $F =\Span \{ \DD_{row}^{1/2} \w_1 , \dots ,\DD_{row}^{1/2} \w_k \}$ and
 $G =\Span \{ \DD_{col}^{1/2} \z_1 , \dots ,\DD_{col}^{1/2} \z_k \}$
with the so-called
normalized row partition vectors $\w_1 ,\dots ,\w_k$ of coordinates
$w_{ji} = \frac1{\sqrt{\Vol (R_i )}}$ if $j\in R_i$ and 0,
otherwise,  and column partition vectors
$\z_1 ,\dots ,\z_k$ of coordinates
$z_{ji} = \frac1{\sqrt{\Vol (C_i )}}$ if $j\in C_i$ and 0, 
otherwise $(i=1,\dots ,k)$. 
Note that the vectors $\DD_{row}^{1/2} \w_1 , \dots ,\DD_{row}^{1/2} \w_k$ and 
$\DD_{col}^{1/2} \z_1 , \dots ,\DD_{col}^{1/2} \z_k$ form orthonormal 
systems in $\R^n$ and $\R^m$, respectively (but they are, usually, not 
complete). By~\cite{Bolla}, we can find  orthonormal
systems  $\tv_1 ,\dots ,\tv_k \in F$ and $\tu_1 ,\dots ,\tu_k \in G$ 
such that
$$
  S_k^2 (\X ) \le \sum_{i=1}^k \| \v_i -\tv_i \|^2 \le 2 S_k^2 (\X ),  \quad
  S_k^2 (\Y ) \le \sum_{i=1}^k \| \uu_i -\tu_i \|^2 \le 2 S_k^2 (\Y ) .
$$
We approximate 
the matrix $\CC_{corr} =\sum_{i=1}^r s_i \v_i \uu_i^T$
by the rank $k$ matrix  $\sum_{i=1}^k s_i \tv_i \tu_i^T$ with the following
accuracy (in spectral norm):
\begin{equation}\label{becsles}
 \left\| \sum_{i=1}^r s_i \v_i \uu_i^T -\sum_{i=1}^k s_i \tv_i \tu_i^T 
 \right\|
 \le \sum_{i=1}^k s_i  \left\| \v_i \uu_i^T -\tv_i \tu_i^T\right\|
 +
 \left\| \sum_{i=k+1}^{r} s_i \v_i \uu_i^T \right\| ,
\end{equation}
where the spectral norm of the last term is at most $\eps$, and the
the individual terms of the first one are estimated from above in the 
following way.
$$
\aligned
 s_i \| \v_i \uu_i^T -\tv_i \tu_i^T \| &\le\| (\v_i  \uu_i^T -\tv_i \uu_i^T ) +
 (\tv_i \uu_i^T -\tv_i \tu_i^T ) \|  \\
&\le\| (\v_i -\tv_i ) \uu_i^T  \| +\| \tv_i (\uu_i -\tu_i )^T  \|  \\
&=\sqrt{\| (\v_i -\tv_i ) \uu_i^T \uu_i (\v_i -\tv_i )^T \| } +
  \sqrt{\|(\uu_i -\tu_i ) \tv_i^T \tv_i (\uu_i -\tu_i )^T \| } \\ 
&=\sqrt{ (\v_i -\tv_i )^T (\v_i -\tv_i  ) } +
  \sqrt{ (\uu_i -\tu_i )^T  (\uu_i -\tu_i ) }  \\
&=\| \v_i -\tv_i  \|  + \| \uu_i -\tu_i  \|  ,
\endaligned
$$
where we exploited that the spectral norm (i.e., the largest singular value)
of an $n\times m$ matrix $\A$ is equal to either the squareroot of the 
largest eigenvalue of the matrix $\A \A^T$ or equivalently, that of $\A^T \A$.
In the above calculations all of these matrices are of rank 1, hence, the
largest eigenvalue of the symmetric, positive semidefinite matrix under the
squareroot is the only non-zero eigenvalue of it, therefore, it is equal to its
trace; finally, we used the commutativity of the trace, and in the last
line we have the usual vector norm.
 
Therefore the first term in~(\ref{becsles}) can be estimated from above by
$$
\aligned
 \sum_{i=1}^k \| \v_i \uu_i^T -\tv_i \tu_i^T \| 
  &\le \sqrt{k} \sqrt{ \sum_{i=1}^k \| \v_i -\tv_i \|^2 } + 
      \sqrt{k} \sqrt{ \sum_{i=1}^k \| \uu_i -\tu_i \|^2 } \\
 &\le \sqrt{k} (\sqrt{2 S_k^2 (\X )} + \sqrt{2 S_k^2 (\Y )})
 = \sqrt{2k} (S_k (\X )+S_k (\Y )). 
\endaligned
$$

Based on these considerations and relation between the cut norm and the
spectral norm (see Lemma~\ref{lemmam}),
the densities to be estimated in the defining formula~(\ref{jeles}) of
volume regularity can be written in terms of stepwise constant vectors in
the following way. The vectors $\hv_i := \DD_{row}^{-1/2} \tv_i$ are stepwise 
constants on the partition $(R_1 ,\dots ,R_k )$ of the rows, whereas 
the vectors $\hu_i := \DD_{col}^{-1/2} \tu_i$ are stepwise 
constants on the partition $(C_1 ,\dots ,C_k )$ of the columns,
$i=1,\dots ,k$. The matrix 
$$
  \sum_{i=1}^k s_i \hv_i \hu_i^T 
$$
is therefore an $n\times m$ block-matrix on $k\times k$
blocks belonging to the above partition of the rows and columns. 
Let ${\hat c}_{ab}$
denote its entries in the $a,b$ block $(a,b=1,\dots ,k)$. 
Using~(\ref{becsles}), the rank $k$ approximation of the matrix $\CC$ is
performed with the following  accuracy of the perturbation $\EE$
in spectral norm:
$$
 \left\| \EE \right\| = \left\| \CC - \DD_{row} (\sum_{i=1}^k 
 s_i \hv_i \hu_i^T ) \DD_{col}  \right\| =
 \left\| \DD_{row}^{1/2} ( \CC_{corr} -\sum_{i=1}^k s_i 
 \v_i \uu_i^T ) \DD_{col}^{1/2} \right\| .
$$
Therefore, the entries of $\CC$ -- for $i\in R_a$, $j\in C_b$ -- 
can be decomposed as
$$
  c_{ij} = d_{row,i} d_{col,j} {\hat c}_{ab} +\eta_{ij} , 
$$
where the cut norm  of 
the $n\times m$ error matrix $\EE =(\eta_{ij} )$  restricted to
$R_a \times C_b$ (otherwise it contains entries all zeroes) 
and denoted by $\EE_{ab}$, is estimated as follows:
$$
\aligned
\| \EE_{ab} \|_{\square} &\le \sqrt{mn} \| \EE_{ab} \| \le \sqrt{nm}\cdot 
 \| \DD_{row,a}^{1/2} \| \cdot 
(\sqrt{2k} (S_k (\X )+S_k (\Y ))+\eps ) \cdot \| \DD_{col,b}^{1/2} \| \\
 &\le \sqrt{nm} \sqrt{c_1 \frac{ \Vol (R_a )}{|R_a |} } \cdot
     \sqrt{c_2 \frac{ \Vol (C_b )}{|C_b |} }  
  (\sqrt{2k} (S_k (\X )+S_k (\Y )) +\eps ) \\
 &= \sqrt{c_1 c_2} \cdot \sqrt{\frac{n}{|R_a|}} \cdot \sqrt{\frac{m}{|C_b|}} 
  \cdot
 \sqrt{\Vol (R_a )} \sqrt{\Vol (C_b )}  (\sqrt{2k} (S_k (\X )+S_k (\Y )) +\eps ) \\
 &\le \sqrt{\frac{c_1 c_2}{K_1 K_2}} \sqrt{\Vol (R_a )} \sqrt{\Vol (C_b )} 
(\sqrt{2k} s 
 +\eps ) \\
 &=c\sqrt{\Vol (R_a )} \sqrt{\Vol (C_b )}  (\sqrt{2k} (S_k (\X )+S_k (\Y )) +\eps ),
\endaligned
$$
where the $n\times n$ diagonal matrix $\DD_{row,a}$ inherits $\DD_{row}$'s
diagonal entries over $R_a$, whereas the $m\times m$ diagonal matrix 
$\DD_{col,b}$ inherits $\DD_{col}$'s diagonal entries over $C_b$,
otherwise they are zeros. Further,
the constants $c_1 ,c_2$ are due to the fact that there are no dominant 
rows and columns, while $K_1 ,K_2$ are derived from the cluster size 
balancing conditions.
Hence,
the constant  $c$ does not depend on $n$ and $m$.
Consequently, for $a,b=1,\dots ,k$ and $X\subset R_a$, $Y\subset C_b$:
$$
\aligned
 &\left| c (X, Y) -\rho (R_a ,C_b ) \Vol (X) \Vol (Y)\right| =\\
 & \left| \sum_{i\in X} \sum_{j \in Y} (d_{row,i} d_{col,j} {\hat c}_{ab}+
 \eta^{ab}_{ij}) - 
\frac{\Vol (X) \Vol (Y )}{\Vol (R_a)  \Vol (C_b )}
\sum_{i\in R_a } \sum_{j \in C_b } (d_{row,i} d_{col,j} {\hat c}_{ab} +
\eta^{ab}_{ij} )   \right| =\\
 &\left| \sum_{i\in X} \sum_{j \in Y} \eta^{ab}_{ij} - 
 \frac{\Vol (X) \Vol (Y )}{\Vol (R_a)  \Vol (C_b )}
 \sum_{i\in R_a } \sum_{j \in C_b } \eta^{ab}_{ij} \right| \le 
 2 \| \EE_{ab} \|_{\square} \\
 &\le 2c(\sqrt{2k} (S_k (\X )+S_k (\Y )) +\eps )\sqrt{\Vol (R_a ) \Vol (C_b )} ,
\endaligned
$$
that gives the required statement for $a,b=1,\dots ,k$.
\end{pf}

Note that when
we use  Definition~\ref{volregcont} of $\gamma$-volume regularity for
the row--column cluster pairs $R_i ,C_j$ $(i,j=1,\dots ,k)$, then we may say
that the \textit{k-way discrepancy} of the underlying contingency table is
the minimum $\gamma$ for which all the row--column cluster pairs are
$\gamma$-volume regular. With this nomenclature, Theorem~\ref{tetelem}
states that the $k$-way discrepancy of a contingency table can be estimated
from above by the  the $(k+1)$th largest singular value of the correspondence
matrix and the $k$-variance of the clusters
obtained by the left and right singular vectors corresponding to the $k$ 
largest singular values of this matrix.
Hence, SVD based representation is
applicable to find volume regular cluster pairs for given $k$, 
where $k$ is the number of structural (protruding) singular values.

\section{Discussion, application, and extension to directed graphs}\label{discussion}

In the ideal $k$-cluster case, we consider the following
generalized random binary contingency table model:
given the partition $(R_1 ,\dots ,R_k )$ of the rows and $(C_1 ,\dots ,C_k )$
of the columns, 
the entry in the row $i\in R_a$ and column  $j\in C_b$ is 1 with
probability $p_{ab}$, and 0 otherwise, independently of other rows of $R_a$ and
columns of $C_b$, $1\le a,b\le k$.
We can think of the probability $p_{ab}$ as the inter-cluster density
of the row--column cluster pair $R_a ,C_b $.
Since generalized contingency tables can be viewed as block-matrices (with
$k\times k$ blocks) burdened with a general random noise, in~\cite{Bol6},
we gave the following spectral characterization of them.
Fixing $k$, and tending with $n$ and $m$ to infinity in
such a way that the cluster sizes grow at the same rate and also $n$ and $m$
subpolynomially, there exists a positive  number $\theta \le 1$, independent
of $n$ and $m$, such that for every $0<\tau <1/2$
there are exactly $k$  singular values of 
$\CC_{corr} $ greater than $\theta -\max \{ n^{-\tau},m^{-\tau} \}$,
while all the others are at most $\max \{ n^{-\tau},m^{-\tau} \}$; further, the
weighted $k$-variance of the row and column representatives constructed 
by the $k$ transformed structural left and right singular vectors is 
${\cal O}(\max \{ n^{-\tau},m^{-\tau} \} )$, respectively.

For general contingency tables,   
our result is that the existence of $k$ singular values of $\CC_{corr}$, 
separated from 0 by $\eps$, is
indication of a $k$-cluster structure, while
the eigenvalues accumulating around 0 are responsible for the
pairwise regularities. The clusters themselves can be recovered by applying
the $k$-means algorithm for the row and column representatives obtained by the
left and right singular vectors corresponding to the structural singular 
values. 

\begin{figure}[ht] \centering
  \includegraphics[scale=0.67]{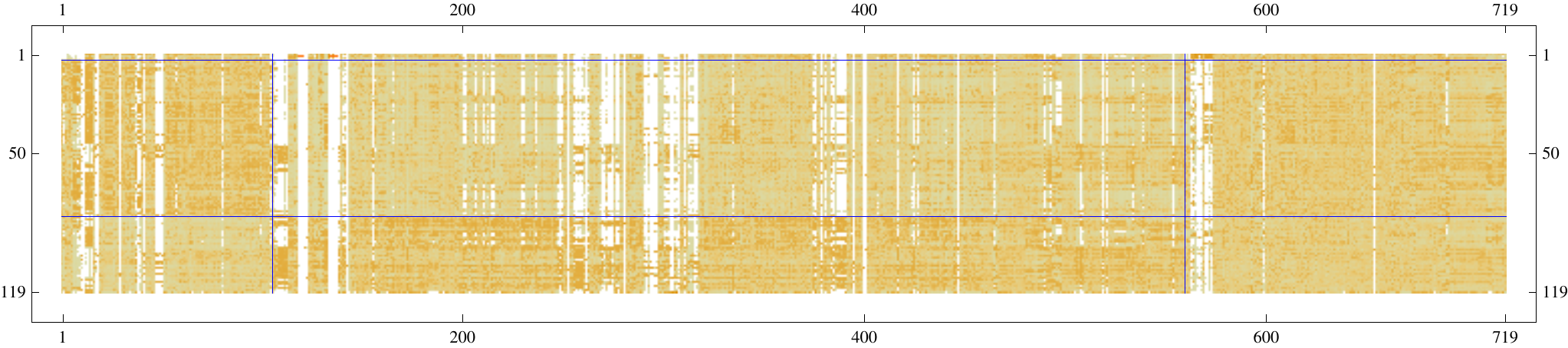}    
\caption[abra]{Result of biclustering 119 stores and 719 products into 3 clusters} 
\label{abra}
\end{figure}

We applied the biclustering algorithm to find simultaneously clusters of
stores and products based on their consumption in TESCO stores. 
Figure~\ref{abra} shows 3 clusters of the stores in which the consumption of the
products belonging to the same cluster was homogeneous with consumption-density
$\frac{c (R_a ,C_b )}{\Vol (R_a ) \Vol (C_b )}$ between store-cluster $R_a$ and
product-cluster $C_b$ $(a,b=1,\dots ,3)$. After sorting the rows and columns 
according to their cluster memberships, we plotted the entries
$\frac{c_{ij}}{d_{row,i} d_{col ,j}}$ (there was one exceptional store-cluster
which contained only 3 stores, but the others could be identified with 
groups of smaller and larger stores associated with product groups of high 
consumption-density within them). 

We can consider quadratic, but
not symmetric contingency tables with zero diagonal as edge-weight matrices
of directed graphs. The $n\times n$ edge-weight matrix $\W$ of a directed
graph has zero diagonal, but is usually not symmetric: $w_{ij}$ is the
weight of the $i\to j$ edge $(i,j=1,\dots ,n; \, i\ne j)$.  
In this setup, the generalized in- and out-degrees are 
$$
 d_{out,i} =\sum_{j=1}^n w_{ij} \quad (i=1,\dots ,n) \quad \textrm{and}\quad
 d_{in,j} =\sum_{i=1}^n w_{ij} \quad (j=1,\dots ,n) ;
$$
further, $\DD_{in} =\diag (d_{in,1}, \dots ,\d_{in,n} )$ and
$\DD_{out} =\diag (d_{out,1}, \dots ,\d_{out,n} )$ are the in- and out-degree
matrices. Suppose that there are no sources and sinks (i.e. no zero out- and
in-degrees), further, that $\W$ is non-decomposable. Then the
correspondence matrix  belonging to $\W$ is
$$
  \W_{corr} = \DD_{out}^{-1/2} \W \DD_{in}^{-1/2} ,
$$
and its SVD is used to minimize the normalized two-way cut of $\W$ as a
contingency table, see Section~\ref{norm}.
Butler~\cite{Butler} generalized the Expander Mixing Lemma for this situation.
We can further generalize it to obtain
regular in- and out-vertex cluster pairs, for a given $k$,
in the following sense.
The $V_{in},V_{out}$ in- and out-vertex cluster pair of the 
directed graph (with sum of the weights of directed edges 1) 
is $\gamma$-volume regular
if for all $X\subset V_{out}$ and $Y\subset V_{in}$ the relation
$$
| w (X, Y) -\rho (V_{out},V_{in}) \Vol_{out} (X) \Vol_{in} (Y)| \le \gamma
 \sqrt{\Vol_{out} (V_{out}) \Vol_{in} (V_{in})} 
$$
holds, where the \textit{directed cut} $w(X,Y)$ is the sum the weights of the 
$X\to Y$ edges,
$\Vol_{out} (X) =\sum_{i\in X} d_{out,i}$, 
$\Vol_{in} (Y) =\sum_{j\in Y} d_{in,j}$, and
$\rho (V_{out},V_{in}) =\frac{w(V_{out},V_{in})}{ \Vol_{out} (V_{out}) \Vol_{in} (V_{in})}$ 
is the relative inter-cluster density of the out--in cluster pair 
$V_{out},V_{in}$.
The clustering $(V_{in,1},\dots ,V_{in,k} )$ and
$(V_{out,1},\dots ,V_{out,k} )$ of the columns and rows -- guaranteed by 
Theorem~\ref{tetelem} -- corresponds to in- and out-clusters of the same
vertex set such that the directed information flow $V_{out,a} \to V_{in,b}$
is as homogeneous as possible for all $a,b=1,\dots ,k$ pairs.

\begin{ack}
We are indebted to the Tesco Hungary for making their data available and
Tam\'as K\'oi for computer processing the data.
\end{ack}


\begin{thebibliography}{00}



\bibitem{Bhatia} Bhatia, R., Matrix Analysis, Springer (1996).

\bibitem{Bilu}
Bilu, Y. and Linial, N., Lifts, discrepancy and nearly optimal 
spectral gap,
\textit{Combinatorica} \textbf{26} (2006), 495--519. 

\bibitem{Bolla} Bolla, M., Tusn\'ady, G., Spectra and optimal partitions of
weighted graphs, \textit{Discrete Mathematics} \textbf{128} (1994), 1--20. 
  


\bibitem{Bol6} Bolla, M., Friedl, K., Kr\'amli, A., Singular value decomposition
of large random matrices (for two-way classification of microarrays),
\textit{Journal of Multivariate Analysis} \textbf{101} (2010), 434--446.



\bibitem{Bolla11} 
Bolla, M., Spectra and structure of weighted graphs,
\textit{Electronic Notes in Discrete Mathematics} \textbf{38} (2011), 149--154.

\bibitem{Bollobas}
Bollob\'as, B., Nikiforov, V., Hermitian matrices and graphs: singular
values and discrepancy, \textit{Discrete Mathematics} \textbf{285} (2004), 
17--32.

\bibitem{Butler}
Butler, S., Using discrepancy to control singular values for nonnegative
matrices, \textit{Lin. Alg. Appl.} \textbf{419} (2006), 486--493.  







\bibitem{Dhillon} Dhillon, I. S., Co-clustering documents and words using
bipartite spectral graph partitioning. In: Proc. ACM Int'l Conf. Knowledge Disc.
Data Mining (KDD 2001), 2001.

\bibitem{Ding} Ding, C., He, X., Zha, H., Gu, M., Simon, H. D.,
A minmax cut spectral method for data clustering and data partitioning,
Lawrence Berkeley National Laboratory Tech. Rep. 54111, 2003.




\bibitem{FKV} Frieze, A., Kannan, R., Vempala, S., Fast Monte-Carlo Algorithms
for finding low-rank approximations. In: Proc. of the 39th Annual IEEE
Symposium on Foundations of Computer Science (FOCS), pp. 370--386, 1998.

\bibitem{Kleinberg} Kleinberg, J., Authoritative sources in hyperlinked
environment, IBM Research Report RJ 10076 (91892), 1997.



\bibitem{Kluger} Kluger, Y., Basri, R., Chang, J. T., Gerstein, M., Spectral
biclustering of microarray data: clustering genes and conditions, 
\textit{Genome Research} \textbf{13} (2003), 703-716. 




\bibitem{Rao} Rao, C. R., Separation theorems for singular values of matrices
and their applications in multivariate analysis, 
\textit{J. Multivariate Analysis} \textbf{9} (1979), 362--377.

\end{thebibliography}
\end{document}